\documentclass[12pt]{article}
\usepackage{oldlfont}
\usepackage{enumerate}
\usepackage{amsmath}
\usepackage{amssymb}
\usepackage{amsthm}
\usepackage{mathrsfs}
\usepackage{amscd}
\usepackage{exscale}
\usepackage{latexsym}
\usepackage[all]{xy}
\usepackage{hyperref}
\usepackage{fancyhdr}
\usepackage{layout}
\usepackage{color}
\usepackage{graphicx}
\DeclareFontFamily{U}{wncy}{}
\DeclareFontShape{U}{wncy}{m}{n}{<->wncyr10}{}
\DeclareSymbolFont{mcy}{U}{wncy}{m}{n}
\DeclareMathSymbol{\Sh}{\mathord}{mcy}{"58} 

\begin{document}

\baselineskip=17pt

\pagestyle{headings}

\numberwithin{equation}{section}

\makeatletter                                                           

\def\section{\@startsection {section}{1}{\z@}{-5.5ex plus -.5ex         
minus -.2ex}{1ex plus .2ex}{\large \bf}}                                 


\pagestyle{fancy}
\renewcommand{\sectionmark}[1]{\markboth{ #1}{ #1}}
\renewcommand{\subsectionmark}[1]{\markright{ #1}}
\fancyhf{} 
\fancyhead[LE,RO]{\slshape\thepage}
\fancyhead[LO]{\slshape\rightmark}
\fancyhead[RE]{\slshape\leftmark}

\addtolength{\headheight}{0.5pt} 
\renewcommand{\headrulewidth}{0pt} 

\newtheorem{thm}{Theorem}
\newtheorem{mainthm}[thm]{Main Theorem}

\newcommand{\ZZ}{{\mathbb Z}}
\newcommand{\GG}{{\mathbb G}}
\newcommand{\Z}{{\mathbb Z}}
\newcommand{\RR}{{\mathbb R}}
\newcommand{\NN}{{\mathbb N}}
\newcommand{\GF}{{\rm GF}}
\newcommand{\QQ}{{\mathbb Q}}
\newcommand{\CC}{{\mathbb C}}
\newcommand{\FF}{{\mathbb F}}

\newtheorem{lem}[thm]{Lemma}
\newtheorem{cor}[thm]{Corollary}
\newtheorem{pro}[thm]{Proposition}
\newtheorem{proprieta}[thm]{Property}
\newcommand{\pf}{\noindent \textbf{Proof.} \ }
\newcommand{\eop}{$_{\Box}$  \relax}
\newtheorem{obs}[thm]{Remark}
\newtheorem{num}{equation}{}

\theoremstyle{definition}
\newtheorem{rem}[thm]{Remark}
\newtheorem*{D}{Definition}

\newcommand{\nsplit}{\cdot}
\newcommand{\G}{{\mathfrak g}}
\newcommand{\GL}{{\rm GL}}
\newcommand{\SL}{{\rm SL}}
\newcommand{\SP}{{\rm Sp}}
\newcommand{\LL}{{\rm L}}
\newcommand{\Ker}{{\rm Ker}}
\newcommand{\la}{\langle}
\newcommand{\ra}{\rangle}
\newcommand{\PSp}{{\rm PSp}}
\newcommand{\U}{{\rm U}}
\newcommand{\GU}{{\rm GU}}
\newcommand{\Aut}{{\rm Aut}}
\newcommand{\Alt}{{\rm Alt}}
\newcommand{\Sym}{{\rm Sym}}

\newcommand{\isom}{{\cong}}
\newcommand{\z}{{\zeta}}
\newcommand{\Gal}{{\rm Gal}}

\newcommand{\F}{{\mathbb F}}
\renewcommand{\O}{{\cal O}} 
\renewcommand{\P}{{\cal P}}
\newcommand{\Q}{{\mathbb Q}}
\newcommand{\R}{{\mathbb R}}
\newcommand{\N}{{\mathbb N}}
\newcommand{\A}{{\mathcal{A}}}
\newcommand{\E}{{\mathcal{E}}}
\newcommand{\J}{{\mathcal{J}}}


\newcommand{\DIM}{{\smallskip\noindent{\bf Proof.}\quad}}
\newcommand{\CVD}{\begin{flushright}$\square$\end{flushright}
\vskip 0.2cm\goodbreak}


\vskip 0.5cm

\title{On the local-global divisibility over ${\rm GL}_2$-type varieties}
\author{Florence Gillibert\footnote{Pontificia Universidad Cat\'olica de Valpara\'so, Supported by Fondecyt Iniciaci\'on Project Number 11130409}, Gabriele Ranieri\footnote{Pontificia Universidad Cat\'olica de Valpara\'so, Supported by Fondecyt Regular Project Number 1140946} }
\date{  }
\maketitle

\vskip 1.5cm

\begin{abstract}
Let $k$ be a number field and let $\A$ be a ${\rm GL}_2$-type variety defined over $k$ of dimension $d$. We show that for every prime number $p$ satisfying certain conditions (see Theorem \ref{teoprinc}), if the local-global divisibility principle by a power of $p$ does not hold for $\A$ over $k$, then there exists a cyclic extension $\widetilde{k}$ of $k$ of degree bounded by a constant depending on $d$ such that $\A$ is $\widetilde{k}$-isogenous to a ${\rm GL}_2$-type variety defined over $\widetilde{k}$ that admits a $\widetilde{k}$-rational point of order $p$. 

Moreover, we explain how our result is related to a question of Cassels on the divisibility of the Tate-Shafarevich group, studied by Ciperiani and Stix and Creutz.
\end{abstract}

\section{Introduction}

Let $k$ be a number field and let ${\mathcal{A}}$ be a commutative algebraic group defined over $ k $.
Several papers have been written on  the following classical question, known as \emph{Local-Global Divisibility Problem}.

\par\bigskip\noindent  P{\small ROBLEM}: \emph{Let $P \in {\mathcal{A}}( k )$. Assume that for all but finitely many valuations $v$ of $k$, there exists $D_v \in {\mathcal{A}}( k_v )$ such that $P = qD_v$, where $q$ is a positive integer. Is it possible to conclude that there exists $D\in {\mathcal{A}}( k )$ such that $P=qD$?}

\par\bigskip\noindent  By  B\'{e}zout's identity, to get answers for a general integer it is sufficient to solve it for powers $p^n$ of a prime. In the classical case of ${\mathcal{A}}={\mathbb{G}}_m$, the answer is positive for $ p $ odd, and negative for instance for $q=8$ (and $P=16$) (see for example \cite{A-T}, \cite{Tro}).

\bigskip For general commutative algebraic groups, Dvornicich and Zannier gave a cohomological interpretation of the problem (see \cite{DZ} and \cite{DZ3}) that we shall explain.
Let $ \Gamma $ be a group and let $ M $ be a $\Gamma$-module.
We say that a cocycle $Z \colon \Gamma \rightarrow M$ satisfies the local conditions if for every $\gamma \in \Gamma$ there exists $m_\gamma \in M$ such that $Z_\gamma = \gamma ( m_\gamma ) - m_\gamma$.
The set of the class of cocycles in $H^1 ( \Gamma, M )$ that satisfy the local conditions is a subgroup of $H^1 ( \Gamma, M )$.
We call it the first local cohomology group $H^1_{{\rm loc}} ( \Gamma, M )$. 
Dvornicich and Zannier \cite[Proposion 2.1]{DZ} proved the following result.

\begin{pro}\label{pro1}
Let $ p $ be a prime number, let $n$ be a positive integer, let $k$ be a number field and let $ \A $ be a commutative algebraic group defined over $ k $.
If $H^1_{{\rm loc}} ( \Gal ( k ( \A[p^n] ) / k ) , \A[p^n] ) = 0$ then the local-global divisibility by $p^n$ over $\A ( k )$ holds.
\end{pro}

The converse of Proposition \ref{pro1} is not true. 
However, in the case when the group $H^1_{{\rm loc}} ( \Gal ( k ( \A[p^n] ) / k ) , \A[p^n] )$ is not trivial, we can find an extension $ L $ of $k$ $k$-linearly disjoint with $k ( \A[p^n] )$ in which the local-global divisibility by $p^n$ over $ \A ( L ) $ does not hold (see \cite[Theorem 3]{DZ3} for the details).

The local-global divisibility problem has been studied for many families of algebraic goups as algebraic tori (see \cite{DZ} and \cite{Ill}), elliptic curves (see \cite{C1}, \cite{C2}, \cite{DZ}, \cite{DZ2}, \cite{DZ3}, \cite{G-R1}, \cite{L-W}, \cite{P1}, \cite{P2}, \cite{P-R-V1}, \cite{P-R-V2}) and, in one of our recent works, polarized abelian surfaces (\cite{G-R2}). 
In particular in \cite{G-R2} we found some algebraic criterions for the triviality of the local cohomology for the torsion points of an abelian variety in the particular case of polarized abelian surfaces. 

In this paper we generalize the algebraic criterions in \cite[Section 2]{G-R2} so that we can apply them to the family of the ${\rm GL}_2$-type varieties, already studied in the particular case of the local-global divisibility problem on torsion points, in \cite{G-R1}.

An abelian variety $ \A $ defined over a number field $ k $ is said to be of ${\rm GL}_2$-type if there exists a number field $ E $ such that $[E: \Q] = {\rm dim} ( \A ) := d$ and an embedding $\phi \colon E \rightarrow {\rm End}_k ( \A ) \otimes \Q$ (see \cite[Section 2, Section 5]{Rib} with a slighty different terminology).
We note ${\mathcal O}_E$ the ring of the algebraic integers of $ E $.
Then $R = \phi^{-1} ( {\rm End}_k ( \A ) \otimes \Z )$ is an order of $\O_E$ and $\phi$ induces an embedding $R \hookrightarrow {\rm End}_k ( \A )$.
Observe that, in certain cases, such field $E$ could be not unique.
We say that a prime number $ p $ is {\bf good} for $ \A $ if it does not divide $[\O_E : R]$.
Suppose also that $ p $ does not ramify over $\O_E$.
For every prime ideal $\P$ dividing $p$, let $\O_{E, \P}$ be the completion of $\O_E$ at $\P$ and consider the representation 
\[
\rho_\P = {\rm Gal} ( \overline{k} / k ) \rightarrow \prod_{\P^\prime \mid p} {\rm GL}_2 ( \O_{E, \P^\prime} ) \rightarrow {\rm GL}_2 ( \O_{E , \P} ),
\]
induced by the action of ${\rm Gal} ( \overline{k} / k )$ over $\A[p^\infty]$.
We prove the following result.

\begin{thm}\label{teoprinc}
Let $k$ be a number field and let $\A$ be a ${\rm GL}_2$-type variety defined over $k$ of dimension $d$. Let $E$ be a field that embeds into ${\rm End}_k ( \A ) \otimes \Q$. Suppose that for every {\bf good} prime number $p$ unramified in $E$ and every prime ideal $\P$ of $\O_E$ over $p$, the determinant of $\rho_\P$ is the cyclotomic character $\chi_p$. There exists an effective constant $C = 2^{( d+1 )^{2 ln ( d ) / ln ( 2 )} + 1}$ depending just on the dimension of $ \A $, such that for every {\bf good} prime number $p \geq 3d + 1$ unramified in $ E $ and such that $k \cap \Q ( \zeta_p ) = \Q$, if there exists $n \in \N$ such that $H^1_{{\rm loc}} ( \Gal ( k ( \A[p^n] ) / k ) , \A[p^n] ) ) \neq 0$, then there exists a ${\rm GL}_2$-type variety $\A^\prime$ isogenous to $\A$ over a finite cyclic extension $\widetilde{k}$ of $k$ such that $[\widetilde{k}: k] \leq C$ and $\A^\prime$ admits a $\widetilde{k}$-rational point of order $p$.
\end{thm}

Let us give some remarks on the hypotheses of Theorem \ref{teoprinc}.
Ribet proved the following result.

\begin{pro}\label{lem23}
Suppose that $\A$ is a ${\rm GL}_2$-type variety defined over a number field $k$ and let $E$ be a totally real field that embeds to ${\rm End}_k ( \A ) \otimes \Q$. Moreover suppose that ${\rm End}_k ( \A ) = {\rm End}_{\overline{k}} ( \A )$. Then ${\rm det} ( \rho_\P ) = \chi_p$, where $\chi_p \colon {\rm Gal} ( \overline{k} / k ) \rightarrow \Z_p^\ast$ is the $p$-adic cyclotomic character.
\end{pro}

\DIM See \cite[Proposition 4.5.1]{Rib}.
\CVD  

Then we have at least a family of ${\rm GL}_2$-type varieties for whom the hypothesis on the determinant of $ \rho_\P $ holds.
Moreover Ribet \cite[p.784]{Rib} observed that Shimura pointed out it would be the case when $E$ is totally real.

Ciperiani and Stix \cite{C-S} and Creutz \cite{C1}, \cite{C2} studied the following question of Cassels which is related to the local-global divisibility problem (see \cite[Remark 20]{C-S} or \cite[Appendix]{G-R2}): let $ k $ be a number field and let $ \A^\prime $ be an abelian variety defined over $ k $. 
For every prime number $q$ we say that the Tate-Shafarevich group $\Sh ( \A^\prime / k )$ is $q$-divisible in $H^1 ( k , \A^\prime )$ if $\Sh ( \A^\prime / k ) \subseteq \cap_{n \in \N^\ast} q^n H^1 ( k , \A^\prime )$. 
What is the set of prime numbers $q$ such that $\Sh ( \A^\prime / k )$ is $q$-divisible ? 

By \cite[Theorem 31]{G-R2} (see also \cite[p. 31]{G-R2}) and Theorem \ref{teoprinc}, we get the following result.

\begin{cor}\label{cor1}
Let $ k $ be a number field and let $\A$ be a principally polarized ${\rm GL}_2$-type variety defined over $k$ of dimension $d$. Let $E$ be a field that embeds into ${\rm End}_k ( \A ) \otimes \Q$. Suppose that for every {\bf good} prime number $p$ unramified in $E$ and every prime ideal $\P$ of $\O_E$ over $p$, the determinant of $\rho_\P$ is the cyclotomic character $\chi_p$. There exists an effective constant $C = 2^{( d+1 )^{2 ln ( d ) / ln ( 2 )} + 1}$ depending on the dimension of $\A$ such that for every {\bf good} prime number $p \geq 3d + 1$ unramified in $E$ and such that $k \cap \Q ( \zeta_p ) = \Q$, if $\Sh ( \A / k )$ is not $p$-divisible in $H^1 ( k , \A )$, then there exists a ${\rm GL}_2$-type variety $\A^\prime$ isogenous to $\A$ over a finite cyclic extension $\widetilde{k}$ of $k$ such that $[\widetilde{k}: k] \leq C$ and $\A^\prime$ admits a $\widetilde{k}$-rational point of order $p$.
\end{cor}

Finally, by using Theorem \ref{teoprinc}, Corollary \ref{cor1} and Weil's bound (see \cite[Theorem 3]{SerreTate} and Section \ref{sec5}), we shall prove the following corollary.

\begin{cor}
\label{BoundIsogeny}
Let $k$ and $E$ be number fields, denote $d=[E:\Q]$. Let $\A/k$ be an abelian variety of dimension $d$ such that $E$ embeds in $ {\rm End}_k(\A) \otimes \Q$. Let $p$ be a prime number. As before, denote $C = 2^{( d+1 )^{2 ln ( d ) / ln ( 2 )} + 1}$. Suppose that : 
\begin{enumerate}
\item  $k \cap \Q ( \zeta_p ) = \Q$;
\item $p$ is a good prime for $\A$;
\item $p\geq 3d+1$;
\item For each place $\P$ of $E$ dividing $p$ the determinant of the representation $\Gal (\overline{k}/k) \rightarrow \GL_{2}(E_{\mathcal{P}})$ is given by the cyclotomic character $\chi_{p}$;
\item $p$ is not ramified in $E$;
\item There exists a place $\mathcal{L}$ of $k$, not dividing $p$ such that $\A$ admits good reduction in $\mathcal{L}$;
\item Denote $\lambda$ the absolute norm of ${\mathcal{L}}$, that is the cardinal of the residue field of $\mathcal{L}$. For each place $\mathcal{P}$ of $E$ dividing $p$, the absolute norm $ {\rm N }(\P)$ of $\P$ satisfies the inequality ${\rm N }(\P) > (1+\lambda^{C/2})^{2d}$. 
\end{enumerate}
Then for all $n \in \N$ we have $H^1_{{\rm loc}} ( \Gal ( k ( \A[p^n] ) / k ) , \A[p^n] ) ) = 0$. Therefore, for all $n \in \N$, the local-global divisibility by $p^n$ holds for $\A ( k )$. If moreover $\A$ is principally polarized, $\Sh ( \A / k )$ is $p$-divisible in $H^1 ( k , \A )$.
\end{cor}
 
Here is the plan of the paper.
In Section \ref{sec2} we recall some well-known results on ${\rm GL}_2$-type varieties and some of our results on local-global divisibility of \cite{G-R1}.

A very important ingredient for the proof of Theorem \ref{teoprinc} is the study of the subgroups of ${\rm GL}_2 ( \F_q )$, where $q$ is a power of a prime number.
We do this in Section \ref{sec3}.

In Section \ref{sec4} we give some algebraic criterions and then we prove Theorem \ref{teoprinc}.

Finally, Section \ref{sec5}, we prove Corollary \ref{BoundIsogeny}.

\section{Definition and basic properties of the ${\rm GL}_2$-type varieties}\label{sec2}

Recall that an abelian variety $\A$ defined over a number field $k$ is said to be of ${\rm GL}_2$-type if there exists a number field $E$ such that $[E: \Q] = {\rm dim} ( \A ) := d$ and an embedding $\phi \colon E \hookrightarrow {\rm End}_k ( \A ) \otimes \Q$ (see \cite[Section 2, Section 5]{Rib}).
We note ${\mathcal O}_E$ the ring of the algebraic integers of $E$.
Then $R = \phi^{-1} ( {\rm End}_k ( \A ) \otimes \Z )$ is an order of $\O_E$ and $\phi$ induces an embedding $R \hookrightarrow {\rm End}_k ( \A )$.
To simplify our notations, we identify $R$ with its image in ${\rm End}_k ( \A ) \otimes \Z$.

We say that a prime number $p$ is {\bf good} for $\A$ if $p$ does not divide the conductor $[\O_E : R]$ (in fact the notion of good depends also on the choice of $E$ if there exists more than a number field holding this property). 
Recall that if $p$ is a good prime number for $\A$, the Tate module $T_p ( \A )$ is a free $\O_E \otimes \Z_p$-module of rank $2$.
Then
\[
\A[p^n] \simeq ( \O_E \otimes \Z_p )^2 / p^n ( \O_E \otimes \Z_p )^2 \simeq \prod_{\P \mid p} ( \O_E / \P^{e_\P n} \O_E )^2,
\]
where $e_\P$ is the ramification index of $\P$.
Each of these direct component of $\A[p^n]$ can also be interpreted in the following way: for every ideal ${\mathcal I}$ of $R$, denote $\A[\mathcal{I}]$ as $\{ x \in \A \ \mid \ \forall \alpha \in {\mathcal I} \ \alpha x = 0 \}$.
Then 
\[
\A[p^n] = \bigoplus_{\P \mid p} \A[\P^{e_\P n}]
\]
with $\A[\P^{e_\P n}] \simeq ( \O_E / \P^{e_\P n} \O_E )^2$.

As a consequence of \cite[Proposition 2.2.1]{Rib}, we have the following result.

\begin{lem}\label{lem21}
Let $p$ be a good prime for $\A$. Then for every positive integer $n$, ${\rm Gal} ( k ( \A[p^n] ) / k )$ is isomorphic to a subgroup of $\prod_{\P \mid p} {\rm GL}_2 ( \O_E / \P^{e_\P n} )$, where $\P$ is a prime of $R$ dividing $p$ and $e_\P$ is the ramification index of $\P$.
\end{lem}

\DIM See \cite[Corollary 15]{G-R1}.
\CVD 

Denote by $\phi_{p, n}$ the embedding 
\[
{\rm Gal} ( k ( \A[p^n] ) / k )  \rightarrow {\rm Aut} ( \A[p^n] ) \simeq \prod_{\P \mid p} {\rm GL}_2 ( \A[\P^{e_\P n}] )
\]
(see also the proof of \cite[Corollary 14]{G-R1}).
Moreover, for every $\P$ prime ideal of $R$ dividing $p$, let $\phi_{\P , n}$ be the composition between $\phi_{p , n}$ and the projection $\prod_{\P^\prime \mid p} {\rm GL}_2 ( \A[\P^{\prime e_{\P} n}] ) \rightarrow {\rm GL}_2 ( \A[\P^{e_\P n}] )$.  
From now on, let us simplify our notation by putting $G_{p, n} = {\rm Gal} ( k ( \A[p^n] ) / k )$ and $G_{\P , n} = \phi_{\P , n} ( G_{p, n} )$. 

The following lemma tells us that, studying the local cohomology of $\A[p^n]$ is equivalent to study the local cohomology of the $G_{\P , n}$-modules $\A[\P^{e_\P n}]$.

\begin{lem}\label{lem22}
Let $p$ be a {\bf good} prime for $\A$ and, for every prime ideal $\P$ of $R$ over $p$, let $e_\P$ be its ramification index. There exists an isomorphism
\[
H^1_{{\rm loc}} ( G_{p, n} , \A[p^n] ) \rightarrow \bigoplus_{\P \mid p} H^1_{{\rm loc}} ( G_{\P , n} , \A[\P^{e_\P n}] ).
\]
\end{lem}

\DIM See \cite[Lemma 16]{G-R1}.
\CVD

Finally, suppose that $p$ is not ramified in $\O_E$ and for every prime ideal $\P$ dividing $p$, let $\O_{E, \P}$ be the completion of $\O_E$ at $\P$. Consider the representation 
\[
\rho_\P = {\rm Gal} ( \overline{k} / k ) \rightarrow \prod_{\P^\prime \mid p} {\rm GL}_2 ( \O_{E, \P^\prime} ) \rightarrow {\rm GL}_2 ( \O_{E , \P} ),
\]
induced by the action of ${\rm Gal} ( \overline{k} / k )$ over $\A[p^\infty]$.
From now on we suppose that the determinant of $\rho_\P$ is the cyclotomic $p$th character $\chi_p$.

\section{Subgroups of ${\rm GL}_2 ( \F_q )$ and reduction to the Borel case}\label{sec3}

By Lemma \ref{lem22} to calculate $H^1_{{\rm loc}} ( G_{p, n} , \A[p^n] )$ is the same that to calculate $H^1_{{\rm loc}} ( G_{\P , n} , \A[\P^{e_\P n}] )$ for any $\P$ dividing $p$.
From now on suppose that $p$ is a {\bf good} prime unramified over $\O_E$.
Then $\A[\P] \simeq ( \O_{E , \P} / \P )^2 \simeq ( \F_q )^2$, where $q$ is a power of $p$ and $q \leq p^d$.
Recall (see Section \ref{sec2}) that the determinant of the representation $\rho_\P$ is the $p$th cyclotomic charcater $\chi_p$. 
In this section we use the classification of the subgroups of ${\rm GL}_2 ( \F_q )$ and some group theory to prove the following proposition.

\begin{pro}\label{pro21}
Suppose that $p \geq 5$, the determinant of the representation $\rho_{\P , 1} \colon G_{\P , 1} \rightarrow {\rm GL}_2 ( \F_q )$ induced by $\rho_\P$ has image equal to $\F_p^\ast$ and that there exists $n \in \N$ such that $H^1_{{\rm loc}} ( G_{\P , n} , \A[\P^n] ) \neq 0$. Then $G_{\P , 1}$ is contained in a Borel subgroup and it is generated by its unique $p$-Sylow subgroup and an element $g$ of order dividing $q-1$, whose determinant has order $p-1$.
\end{pro}

\DIM To prove Proposition \ref{pro21}, we prove some several preliminary lemmas.
First, let us prove that $G_{\P, 1}$ is isomorphic to its projective image.

\begin{lem}\label{lem222}
$G_{\P, 1}$ is isomorphic to its projective image in ${\rm PGL}_2 ( \F_q )$.
\end{lem}

\DIM Suppose that it is not the case. Then $G_{\P , 1}$ contains an element that is a non-trivial scalar multiple of the identity (i.e. an element which is in any basis diagonal with the same eigenvalues distinct from $1$).
Then, use the argument in \cite[p. 29]{DZ3} to prove that $H^1 ( G_{\P , n} , \A[\P^n] ) = 0$.
\CVD

We now study the case when $p$ does not divide the order of $G_{\P, 1}$.

\begin{lem}\label{lem233}
Suppose that $p$ does not divide the order of $G_{\P , 1}$. Then $G_{\P , 1}$ is cyclic, generated by an element $g$ of order dividing $q-1$, whose determinant has order $p-1$. In particular it is contained in a Borel subgroup.
\end{lem}   

\DIM By Lemma \ref{lem222}, $G_{\P, 1}$ is isomorphic to its projective image. 
It is well-kwown, see for instance \cite[Proposition 16]{Ser}, that all the possible subgroups of ${\rm PGL}_2 ( \F_q )$ of order not divided by $p$ are either cyclic, or dihedral, or isomorphic to one of the following three groups: $A_4$, $S_4$, $A_5$.

All dihedral group may be generated by elements of order $2$.
Then for $p \geq 5$ this contradicts the fact that the image of the determinant of $\rho_{\P , 1}$ is $\F_p^\ast$.

If $G_{\P , 1}$ is isomorphic either to $A_4$, or to $S_4$ or to $A_5$, then we get a contradiction because in this case $G_{\P , 1}$ contains a subgroup isomorphic to $\Z / 2 \Z \times \Z / 2 \Z$.
Since $G_{\P , 1}$ is a subgroup of ${\rm GL}_2 ( \F_q )$, we get that $-Id \in G_{\P , 1}$ and this is not possible because $G_{\P , 1}$ is isomorphic to its projective image.

Then $G_{\P , 1}$ is cyclic and let $g$ be its generator.
Suppose that $g$ has order not dividing $q-1$.
Since $p$ does not divide the order of $G_{\P , 1}$,
the eigenvalues of $g$ are in $\F_{q^2}$ and there are of the form $\lambda$, $\lambda^q$, with $\lambda^{q+1} \in ( \Z / p \Z )^\ast$.
Since the determinant of $\rho_{\P , 1}$ has image $\F_p^\ast$, then $\lambda^{q+1}$ is distinct from $1$ and so $g^{q+1}$ is a scalar matrix distinct from the identity. 
This is not possible by Lemma \ref{lem222}.
\CVD 

We now study the case when $p$ divides the order of $G_{\P, 1}$.

\begin{lem}\label{lem244}
Suppose that $p$ divides the order of $G_{\P, 1}$. Then $G_{\P, 1}$ has an elementary abelian normal $p$-Sylow subgroup $N$. Moreover $G_{\P, 1}$ is generated by $N$ and an element of order dividing $q-1$ and determinant of order $p-1$. In particular, $G_{\P, 1}$ is contained in a Borel subgroup.
\end{lem}

\DIM Suppose that $G_{\P, 1}$ has more than one $p$-Sylow subgroup.  
Then by the classification of the subgroups of ${\rm SL}_2 ( \F_q )$, see for instance \cite[Chapter 3, Theorem 6.17]{Suz}, $G_{p, 1}$ contains ${\rm SL}_2 ( \F_p )$. 
By Lemma \ref{lem222} this is not possible.
Then $G_{\P, 1}$, see again \cite[Chapter 3, Theorem 6.17]{Suz}, contains a unique $p$-Sylow subgroup $N$, which is elementary abelian.
Since $N$ is normal, $G_{\P , 1}$ is contained in a Borel subgroup.
Let $\tau$ be an element of $G_{\P , 1}$.
As $G_{\P , 1}$ is in a Borel subgroup, the eigenvalues of $\tau$ are in $\F_q$. 
If $\tau$ has two distinct eigenvalues, it is diagonalizable and $\tau^{q-1} = Id$.
If $\tau$ has a unique eigenvalue $\alpha$, then $\tau^p = \alpha Id$ and, by Lemma \ref{lem222}, this is possible only if $\alpha = 1$.

Since the determinant of $\rho_{\P , 1}$ has image $\F_p^\ast$, $G_{\P , 1}$ has an element $g$ of order dividing $q-1$ and with determinant of order $p-1$.
Suppose that $g$ and $N$ does not generate $G_{\P, 1}$.
Then there exists $g^\prime$ of order dividing $q-1$ such that $g^\prime N$ is not a power of $g N$.
Since $G_{\P, 1} / N$ is isomorphic to a subgroup of $( \F_q^\ast )^2$, this is possible only if for a certain prime $l$ dividing $q-1$, $G_{\P, 1}$ contains a subgroup isomorphic to $( \Z / l \Z )^2$. Then $G_{\P, 1}$ contains a non-trivial scalar matrix and this contradicts Lemma \ref{lem222}.
\CVD

Proposition \ref{pro21} immediately follows from Lemma \ref{lem233} and Lemma \ref{lem244}.
\CVD

\section{Some algebraic criterions and proof of Theorem \ref{teoprinc}}\label{sec4}

In \cite[section 2]{G-R2} we proved some criterions for the triviality of the local cohomology for a large class of Galois-module. 
In this section we generalize these criterions, and we apply them on the family of the ${\rm GL}_2$-type varieties.

From now on we suppose that there exists $n \in \N$ such that the cohomology group $H^1_{{\rm loc}} ( G_{\P , n} , \A[\P^n] )$ is distinct from $0$. 
Then, by Proposition \ref{pro21}, $G_{\P, 1}$ is contained in a Borel subgroup. 
Recall that there exists $q$ a power of $p$ less or equal than $p^d$ such that $\A[\P]$ is isomorphic to $( \F_q^\ast )^2$. 
Then, again by Proposition \ref{pro21}, $G_{\P, 1}$ is generated by its unique $p$-Sylow subgroup $N$, which is either trivial or elementary abelian of exponent $p$, and by an element $g$ of order dividing $q-1$ and whose determinant has order $p-1$.
Let us denote by $g_n$ an element of $G_{\P , n}$ whose image by the projection $\pi_n \colon G_{\P, n} \rightarrow G_{\P, 1}$ is $g$.
We can suppose that the order of $g_n$ is equal to the order of $g$ by replacing $g_n$ by a suitable $p$-power.

The following lemma gives some conditions on the eigenvalues of $g$.

\begin{lem}\label{lem31}
Let $i$ be ${\rm ord} ( g )/ ( p - 1 )$. Then $g^i$ has an eigenvalue equal to $1$.
\end{lem}

\DIM Consider the projection $\pi_n \colon G_{\P , n} \rightarrow G_{\P , 1}$.
Let $\widetilde{N}$ be $\pi_n^{-1} ( N )$.
Since $\ker ( \pi_n )$ is a $p$-group and $N$ is the $p$-Sylow of $G_{\P, 1}$, $\widetilde{N}$ is the $p$-Sylow of $G_{\P , n}$.
Suppose that $g^i$ has the eigenvalues distinct from $1$ and let $g_n$ be as before.
Since $g = \pi_n ( g_n )$ and $g_n$ has the same order of $g$, then $g_n^i - Id \colon \A[\P^n] \rightarrow \A[\P^n]$ is bijective.   
Now, just copy the proof of \cite[Corollary 9]{G-R2} with $g_n$ in the place of $g$ and $\widetilde{N}$ in the place of $H$.
\CVD
 
\begin{rem}\label{rem32}
Since it will be very useful, let us obtain a bound for $( i , p-1 )$. Recall that $g$ has order which divides $q-1$ and $q$ is a power of $p$ equals or less than $p^d$, where $d$ is the dimension of $\A$. Suppose that $q = p^b$. Then an easy calculation shows that $( ( p^b - 1 ) / ( p-1 ) , p-1 ) \leq b$. Then $( i, p-1 ) \leq d$. 
In particular the bound does not depend on $p$, but only on the dimension of the variety.
\end{rem}

Now we want to get some other restrictions to the eigenvalues of $g$, by using \cite[Lemma 13, Proposition 17]{G-R2}.

\begin{lem}\label{lem33}
Suppose that the eigenvalues of $g$ are distinct from $1$. Then the homomorphism $H^1 ( G_{\P , n} , \A[\P] ) \rightarrow H^1 ( G_{\P , n} , \A[\P^n] )$ induced by the exact sequence of $G_{\P , n}$-modules
\[
0 \rightarrow \A[\P] \rightarrow \A[\P^n] \rightarrow \A[\P^{n-1}] \rightarrow 0
\]
is injective and its image is $H^1 ( G_{\P , n} , \A[\P^n] )[p]$. In other words $H^1 ( G_{\P , n} , \A[\P] )$ is isomorphic to $H^1 ( G_{\P , n} , \A[\P^n] )[p]$.
\end{lem}

\DIM Observe that since $p$ is a {\bf good} prime unramified over $E$, we have $\A[p] \cap \A[\P^{n-1}] = \A[\P]$. 
Then the following sequence 
\[ 
0 \rightarrow \A[\P] \rightarrow \A[\P^n] \rightarrow \A[\P^{n-1}] \rightarrow 0 
\]
(here the first map is the inclusion and the second is the multiplication by $p$) is exact.
Copy the proof of \cite[Lemma 13]{G-R2} with $G_{\P, n}$ in the place of $G$, $\A[\P]$, $\A[\P^n]$, $\A[\P^{n-1}]$ in the place respectively of $V_{n, d}[p]$, $V_{n, d}$, $V_{n, d}[p^{n-1}]$ and $g_n$ in the place of $\delta$.
\CVD

\begin{lem}\label{lem34}
Suppose that the eigenvalues of $g$ are distinct from $1$. If $H^1 ( G_{\P, 1} , \A[\P] ) = 0$ and $\A[\P]$ and ${\rm End} ( \A[\P] )$ have no common irreducible sub $\Z / p \Z [\langle g \rangle]$-module, then $H^1 ( G_{\P, n}, \A[\P^n] ) = 0$, where $g$ acts by conjugacy over ${\rm End} ( \A[\P] )$. 
\end{lem}

\DIM Let $H$ be the normal subgroup of $G_{\P, n}$ of the elements acting like the identity over $\A[\P]$. 
Observe that $G_{\P, 1}$ is isomorphic to $G_{\P, n} / H$ and consider the inflation-restriction sequence:
\[
0 \rightarrow H^1 ( G_{\P, 1} , \A[\P] ) \rightarrow H^1 ( G_{\P, n} , \A[\P] ) \rightarrow H^1 ( H , \A[\P] )^{G_{\P, 1}}. 
\]
By copying the proof of \cite[Proposition 17]{G-R2} with $G_{\P, n}$ in the place of $G$, $g_n$ in the place of $\delta$, $g$ in the place of $\overline{\delta}$ and $\A[\P]$ in the place of $V_{n, d}[p]$, we get $H^1 ( G_{\P, n} , \A[\P] ) = 0$.
By Lemma \ref{lem33}, we get $H^1 ( G_{\P , n} , \A[\P^n] )[p] = 0$, and so $H^1 ( G_{\P, n}, \A[\P^n] ) = 0$.
\CVD

In the previous lemma we made the hypothesis that $H^1 ( G_{\P, 1} , \A[\P] ) = 0$.
Next lemma shows that actually this condition is already implied by the other conditions when $p \geq 3d + 1$.

\begin{lem}\label{lem35}
Suppose that the eigenvalues of $g$ are distinct from $1$ and that $H^1 ( G_{\P, 1} , \A[\P] ) \neq 0$. Moreover suppose that $p \geq 3d + 1$. Then $\A[\P]$ and ${\rm End} ( \A[\P] )$ have a common irreducible sub $\Z / p \Z [\langle g \rangle]$-module.
\end{lem}

\DIM Recall (see Proposition \ref{pro21}), that $\A[\P]$ is a $\F_q$-vector space of dimension $2$, $G_{\P , 1}$ is included in a Borel subgroup of ${\rm GL}_2 ( \F_q )$ and it contains a unique $p$-Sylow subgroup that we denote by $N$. 
So $\A[\P]$ has a non-trivial sub $G_{\P , 1}$-module $V$ (which is a sub $\F_q$-vector space of dimension $1$) such that $N$ acts like the identity over $V$ and $\A[\P]/ V$.
The following exact sequence of $G_{\P , 1}$-modules
\[
0 \rightarrow V \rightarrow \A[\P] \rightarrow \A[\P] / V \rightarrow 0,
\]
gives rise to the following long cohomology exact sequence:
\[
\cdots H^0 ( G_{\P, 1} , \A[\P] / V ) \rightarrow H^1 ( G_{\P , 1} , V ) \rightarrow H^1 ( G_{\P, 1} , \A[\P] ) \rightarrow H^1 ( G_{\P, 1} , \A[\P] / V ).    
\]
Since $g$ has all the eigenvalues distinct from $1$, $H^0 ( G_{\P, 1} , \A[\P] / V ) = 0$, and so we have the exact sequence
\begin{equation}\label{rel1}
0 \rightarrow H^1 ( G_{\P , 1} , V ) \rightarrow H^1 ( G_{\P, 1} , \A[\P] ) \rightarrow H^1 ( G_{\P, 1} , \A[\P] / V ).
\end{equation}
Since $N$ is normal in $G_{\P , 1}$, we have the inflation-restriction sequences:
\[
0 \rightarrow H^1 ( G_{\P , 1}/N , V ) \rightarrow H^1 ( G_{\P, 1} , V ) \rightarrow H^1 ( N , V )^{G_{\P , 1} / N}.
\]
and
\[
0 \rightarrow H^1 ( G_{\P , 1}/N , \A[\P] / V ) \rightarrow H^1 ( G_{\P, 1} , \A[\P] / V ) \rightarrow H^1 ( N , \A[\P] / V )^{G_{\P , 1} / N}.
\]
Since $N$ is a $p$-Sylow subgroup of $G_{\P , 1}$, the map $H^1 ( G_{\P, 1} , V ) \rightarrow H^1 ( N , V )^{G_{\P , 1} / N}$ and the map $H^1 ( G_{\P, 1} , \A[\P] / V ) \rightarrow H^1 ( N , \A[\P] / V )^{G_{\P , 1} / N}$ are injective.
Since $N$ acts like the identity over $V$ and $\A[\P] / V$, $N$ and $V$ are abelian $p$-groups of exponent $p$ and $G_{\P , 1} / N$ is generated by the class of $g$ modulo $N$, we have $H^1 ( N , V )^{G_{\P , 1} / N} = {\rm Hom}_{\Z / p \Z[\langle g \rangle]} ( N , V )$, and $H^1 ( N , \A[\P] / V )^{G_{\P , 1} / N} = {\rm Hom}_{\Z / p \Z[\langle g \rangle]} ( N , \A[\P] / V )$.
By Lemma \ref{lem31} and Remark \ref{rem32}, we have that for $p \geq 3d + 1$, there exists $i$ such that $g^i$ has eigenvalues $\lambda$ and $1$ with $\lambda \in \F_p^\ast$ of order $\geq 3$.
Let $h$ be any element of $N$ and let $Z$ be in ${\rm Hom}_{\Z / p \Z[\langle g \rangle]} ( N , \A[\P] / V )$.
Then
\begin{equation}\label{rel2}
Z_{g^i h g^{-i}} = g^i Z_h.
\end{equation}
Suppose first that $g^i$ acts like the identity over $\A[\P] / V$.
Then $g^i h g^{-i} = h^{\lambda}$ and so, by relation (\ref{rel2}), we get $\lambda Z_h = Z_h$.
Thus, since $\lambda \neq 1$, for every $h \in N$ we have $Z_h = 0$.

Now suppose that $g^i$ acts like the multiplication by $\lambda$ over $\A[\P] / V$.
Then $g^i h g^{-i} = h^{\mu}$, where $\mu \in \F_p^\ast$ and $\lambda \mu = 1$. 
So, we get $\mu Z_h = \lambda Z_h$.
Thus, since $\lambda$ has order $\geq 3$, for every $h \in N$, $Z_h = 0$.
Then ${\rm Hom}_{\Z / p \Z[\langle g \rangle]} ( N , \A[\P] / V ) = 0$.
Thus $H^1 ( G_{\P, 1} , \A[\P] / V ) = 0$ and, by relation (\ref{rel1}), we get that $H^1 ( G_{\P , 1} , \A[\P] )$ is isomorphic to $H^1 ( G_{\P, 1} , V )$. 
Since we proved that $H^1 ( G_{\P , 1} , V )$ embeds into ${\rm Hom}_{\Z / p \Z[\langle g \rangle]} ( N , V )$, we have that $H^1 ( G_{\P , 1} , \A[\P] )$ is isomorphic to a subgroup of ${\rm Hom}_{\Z / p \Z[\langle g \rangle]} ( N , V )$.
Since $H^1 ( G_{\P , 1} , \A[\P] ) \neq 0$, we have a non-trivial $Z \colon N \rightarrow V$ homomorphism of $\Z / p \Z[\langle g \rangle]$-modules.
Let $l ( N )$ be the sub-space of ${\rm End} ( \A[\P] )$ given by $\{ h -Id \ \mid \ h \in N \}$.
Then $l ( N )$ is a $\Z / p \Z [\langle g \rangle]$-module and $Z$ induces a non-trivial homomorphism $\phi_Z \colon l ( N ) \rightarrow V$ by sending $h - Id$ to $Z_h$.
Thus a quotient of $l ( N )$ is isomorphic to a sub-module of $V$ (then of $\A[\P]$).
But since the action of $g$ is semisimple, every quotient is isomorphic to a submodule. 
Hence $\A[\P]$ and ${\rm End} ( \A[\P] )$ have a common irreducible sub$\Z / p \Z [\langle g \rangle]$-module.
\CVD

Lemmas \ref{lem34} and \ref{lem35} tell us that, under the hypothesis $p \geq 3d + 1$, $\A[\P]$ and ${\rm End} ( \A[\P] )$ have a common sub $\Z / p \Z[\langle g \rangle]$-module.
We use this in the next corollary to get a restriction on the possible eigenvalues of $g$.

\begin{cor}\label{cor36} 
Let $\lambda_1$, $\lambda_2 \in \F_q^\ast$ be the eigenvalues of $g$ and let $\mu_1$, $\mu_2$ be in $\F_q^\ast$ such that $\lambda_1 \mu_1 = 1$ and $\lambda_2 \mu_2 = 1$. Suppose that $\lambda_1 \neq 1$ and $\lambda_2 \neq 1$ and that $p \geq 3d + 1$. Then at least one of the following holds:
\begin{enumerate}
\item $\lambda_1$ is Galois conjugated to $\lambda_1 \mu_2$;
\item $\lambda_1$ is Galois conjugated to $\lambda_2 \mu_1$;
\item $\lambda_2$ is Galois conjugated to $\lambda_1 \mu_2$;
\item $\lambda_2$ is Galois conjugated to $\lambda_2 \mu_1$.
\end{enumerate}
Moreover, either $\lambda_1$ or $\lambda_2$ has order bounded by $C = 2^{( d+1 )^{2 ln ( d ) / ln ( 2 )} + 1}$.
\end{cor}

\DIM Observe that $\A[\P] \simeq \F_q^2$ is the direct sum of two $\F_q$-vector spaces $V_1$ and $V_2$ of dimension $1$, and, for $i = 1, 2$, $g$ acts on $V_i$ like the multiplication by $\lambda_i$.
In particular $V_1$ and $V_2$ are two $\Z / p \Z[\langle g \rangle]$-modules. 
For $i \in \{1 , 2 \}$, the minimal polynomial of the restriction of $g$ to $V_i$ over $\F_p$ is equal to the minimal polynomial of $\lambda_i$ over $\F_p$.  
On the other hand ${\rm End} ( \A[\P] )$ is the direct sum of three $\Z / p \Z[\langle g \rangle]$-modules 
$T_1$, $T_2$ and $T_3$.
We can suppose that the automorphism of ${\rm End} ( \A[\P] )$ associated to $g$ has minimal polynomial over $T_1$ equal to the minimal polynomial of $\lambda_1 \mu_2$ over $\F_p$, it is the identity over $T_2$ and it has minimal polynomial over $T_3$ equal to the minimal polynomial of $\lambda_2 \mu_1$ over $\F_p$. By Lemmas \ref{lem34} and \ref{lem35}, $\A[\P]$ and ${\rm End} ( \A[\P] )$ have a non-trivial common $\Z / p \Z[\langle g \rangle]$-module.
Then, since $\lambda_1 \neq 1$ and $\lambda_2 \neq 1$, this is possible only if one of the fourth condition of this corollary holds.

Now recall that by Proposition \ref{pro21} $\lambda_1 \lambda_2 \in \F_p^\ast$ and so $\lambda_1^{p-1} \lambda_2^{p-1} = 1$.
Suppose that condition 1. holds (or equivalently condition 4. holds, because condition 4. is identical to condition 1. if we replace $\lambda_1$ with $\lambda_2$).
Then $\lambda_1$ is conjugated to $\lambda_1 \mu_2$.
Thus there exists $l < d$ such that $\lambda_1^{( p^l -1 )} = \mu_2$. 
Now recall that by Proposition \ref{pro21}, $\lambda_1 \lambda_2 \in \F_p^\ast$ and so $\lambda_1^{p-1} \lambda_2^{p-1} = 1$.
Then since $\mu_2$ is the inverse of $\lambda_2$, we get that $\lambda_1^{p^l - 2}$ has order dividing $p-1$.
Again because $\lambda_1 \lambda_2 \in \F_p^\ast$, we have $\lambda_2^{p^l -2}$ has order dividing $p-1$.
Then $g^{p^l - 2}$ has order dividing $p-1$ and so, by Lemma \ref{lem31}, $g^{p^l -2}$ has an eigenvalue equal to $1$.
On the other hand $g$ has order dividing $q-1$ and $q = p^b$ with $b \leq d$.
Then either $\lambda_1$ or $\lambda_2$ has order dividing $( p^b -1 , p^l - 2 )$.

Suppose now that condition 2. holds (or equivalently condition 3. holds, because condition 3. is identical to condition 2. if we replace $\lambda_1$ with $\lambda_2$).
A calculation identical to the previous, gives that there exists $r < d$ such that $\lambda_1^{p^r +2}$ has order dividing $p-1$.
Then as in the previous case we get that $g^{p^r + 2}$ has order dividing $p-1$ and so, by Lemma \ref{lem31}, $g^{p^r + 2}$ has an eigenvalue equal to $1$.
Then either $\lambda_1$ or $\lambda_2$ has order dividing $( p^b - 1 , p^r + 2 )$.

The bound on the order of either $\lambda_1$ or $\lambda_2$ now follows from the following lemma (applied in the case $N = 2$).

\begin{lem}
\label{Euclid}
Let $a\in \N$ and $b\in \N$ be two integers less than $d$, let $A,B,C,D$ be relative integers all with absolute values less than $N \in \N^\ast$. The greatest commun divisor of $Ap^a-B$ and $Cp^b-D$ is less than $ 2 N^{(d+1)^{2 ln(d)/ln(2)}}.$
\end{lem}
\DIM Let us make an induction on the exponents $a$ and $b$. Without loss of generality we suppose $b\le a$, let $q$ and $0\le r<b$ be the quotient and the rest of the Euclidian division of $a$ by $b$ (so $a=bq+r$). Observe that : 
$$ (Ap^a-B,Cp^b-D) \le (C^{q}Ap^{bq+r}-C^{q}B,Cp^b-D), $$
but we have :
$$ C^{q}Ap^{bq+r}-C^{q}B \equiv AD^qp^r-C^{q}B \mod ( Cp^b-D ) ,$$
so we can bound $(Ap^a-B,Cp^b-D)$ by $( AD^qp^r-C^{q}B,  Cp^b-D)$.
Observe that all the coefficients $AD^q$, $C^{q}B$ are less than $N^{d+1}$. As it is well known that the Euclid algorithm applied to $a$ by $b$ will finish in at most $2 ln(d)/ln(2)$ steps, we can bound the final result by $ 2 N^{(d+1)^{2 ln(d)/ln(2)}} $. 
\CVD

\CVD

We can now prove Theorem \ref{teoprinc}.

\begin{thm}[{\rm Theorem \ref{teoprinc}}]
Let $k$ be a number field and let $\A$ be a ${\rm GL}_2$-type variety defined over $k$ of dimension $d$. Let $E$ be a field that embeds into ${\rm End}_k ( \A ) \otimes \Q$. Suppose that for every {\bf good} prime number $p$ unramified in $E$ and every prime ideal $\P$ of $\O_E$ over $p$, the determinant of $\rho_\P$ is the cyclotomic character $\chi_p$. There exists an effective constant $C = 2^{( d+1 )^{2 ln ( d ) / ln ( 2 )} + 1}$ depending just on the dimension of $ \A $, such that for every {\bf good} prime number $p \geq 3d + 1$ unramified in $ E $ and such that $k \cap \Q ( \zeta_p ) = \Q$, if there exists $n \in \N$ such that $H^1_{{\rm loc}} ( \Gal ( k ( \A[p^n] ) / k ) , \A[p^n] ) ) \neq 0$, then there exists a ${\rm GL}_2$-type variety $\A^\prime$ isogenous to $\A$ over a finite cyclic extension $\widetilde{k}$ of $k$ such that $[\widetilde{k}: k] \leq C$ and $\A^\prime$ admits a $\widetilde{k}$-rational point of order $p$.
\end{thm}

\DIM Let $p$ be a prime number that satisfies all the hypotheses of the Theorem.
Since $k \cap \Q ( \zeta_p ) = \Q$ the determinant of the representation $\rho_{\P , 1} \colon G_{\P , 1} \rightarrow {\rm GL}_2 ( \F_q )$ induced by $\rho_\P$ has image equal to $\F_p^\ast$ (see Proposition \ref{pro21}).   
Then, by Proposition \ref{pro21}, there exists a prime ideal $\P$ of ${\cal O}_E$ such that $G_{\P, 1}$ (see Lemmas \ref{lem21} and \ref{lem22}) contains an element $g$ of order dividing $q-1$ (where $q$ is the inertia degree of $\P$) and with determinant of order $p-1$, such that $G_{\P, 1}$ is generated by $g$ and its unique $p$-Sylow subgroup $N$.  
Moreover by Corollary \ref{cor36} there exists an effective constant $c \leq 2^{( d+1 )^{2 ln ( d ) / ln ( 2 )} + 1}$ depending just on $d$ such that $g^c$ has an eigenvalue equal to $1$.
Recall (see p.6 after Lemma \ref{lem21}) that $G_{\P , 1}$ is by definition $\phi_{\P, 1} ( {\rm Gal} ( k ( \A[p] ) / k ) )$, where $\phi_{\P , 1}$ is the composition between the embedding
\[
{\rm Gal} ( k ( \A[p] ) / k )  \rightarrow {\rm Aut} ( \A[p] ) \simeq \prod_{\P \mid p} {\rm GL}_2 ( \A[\P] )
\]
and the projection $\prod_{\P^\prime \mid p} {\rm GL}_2 ( \A[\P^{\prime}] ) \rightarrow {\rm GL}_2 ( \A[\P] )$.  
Then set $\widetilde{k}$ the subfield of $k ( \A[p] )$ fixed by $\phi_{\P , 1}^{-1} ( \langle g^c , N \rangle )$.   
Observe that $\langle g^c , N \rangle$ is a normal subgroup of $G_{\P , 1}$ and that ${\rm Gal} ( k ( \A[p] ) / k ) / \phi_{\P , 1}^{-1} ( \langle g^c , N \rangle )$ is isomorphic to $G_{\P , 1} / \langle g^c , N \rangle$.  
Hence $\widetilde{k}$ is a cyclic extension of $k$ of degree $c$ contained in $k ( \A[p] )$ such that $\Gal ( k ( \A[p] ) / \widetilde{k} )$ acts over $\A[\P]$ as the group generated by $N$ and $g^c$.
Thus $\A[\P]$ has a non-trivial Galois sub-module $V$ in which $N$ acts like the identity and $g^c$ acts either like the identity or by multiplication by an element of $\F_p^\ast$.
If $g^c$ acts like the identity over $V$ then every non-trivial element of $V$ is a $\widetilde{k}$-rational point of order $p$ of $\A$.

If $g^c$ does not act like the identity, let $\A^\prime$ be the ${\rm GL}_2$-type variety defined over $\widetilde{k}$ and $\widetilde{k}$-isogenous to $\A$ with an isogeny with kernel $V$.
Then, since $g^c$ and $N$ act like the identity over $\A[\P] / V$, $\A^\prime$ admits a $\widetilde{k}$-rational point of order $p$.
\CVD

\section{Proof of Corollary \ref{BoundIsogeny}}\label{sec5}

The proof of Corollary \ref{BoundIsogeny} is a consequence of Theorem \ref{teoprinc} and of Weil's results on the eigenvalues of the images of Frobenius substitutions for primes of good reduction, by the Galois representation associated to $\A$ (cf. Theorem \ref{WeilBound} below). In fact Weil's results also worth for the $\GL_2$-type representation.  

Let $\overline{{\mathcal{L}}}$ be a place of $\overline{k}$ over ${\mathcal{L}}$. We denote $D(\overline{\mathcal{L}})=\{  \sigma \in\Gal(\overline{k}/k) \ : \ \sigma(\overline{\mathcal{L}})=\overline{\mathcal{L}}  \} $ and 
$I(\overline{\mathcal{L}})=\{  \sigma \in D(\overline{\mathcal{L}}) : \  \forall x\in \overline{k} \ \sigma(x) \equiv x \mod \overline{\mathcal{L}}  \} $ the decomposition group and the inertia group respectively. We have a natural morphism $D(\overline{\mathcal{L}}) \rightarrow \Gal ( \overline{\F}_{\lambda} / \F_{\lambda}) $ whose kernel is $I(\overline{\mathcal{L}})$. By the results of Weil (see for instance \cite[Theorem 3]{SerreTate}) we have the following theorem:
\begin{thm}
\label{WeilBound}
Suppose that $\mathcal{A}$ has potentially good reduction at the prime $\mathcal{L}$ of $k$. Consider the representation associated to the Tate module $T_p(A)$ : 
$$\rho_p \colon \Gal(\overline{k}/k) \rightarrow \GL_{2d}(\Z_p).$$
Let $\sigma$ be an element of $D(\overline{\mathcal{L}})$ such that the reduction of $\sigma$ in $\Gal ( \overline{\F}_{\lambda} / \F_{\lambda}) $ is the Frobenius morphism. Then the characteristic polynomial of $\rho_p(\sigma)$ has integral coefficients independent of $p$. The absolute values of its roots are equal to ${\rm N}({\mathcal{L}})^{1/2}=\lambda^{1/2}$. 
\end{thm}
 
This also worth for the $\GL_2$-type representation.    
\begin{cor}
\label{Weil}
Suppose that $\mathcal{A}$ has good reduction at the prime $\mathcal{L}$ of $k$. Consider the representation of $\GL_2$-type  
$$\rho_{\mathcal{P}} \colon \Gal(\overline{k}/k ) \rightarrow \GL_2(O_{E,\mathcal{P}}).$$
Let $\sigma$ be an element of $D(\overline{\mathcal{L}})$ such the reduction of $\sigma$ in $\Gal ( \overline{\F}_{\lambda} / \F_{\lambda}) $ is the Frobenius morphism. Then the characteristic polynomial of $\rho_{\mathcal{P}} (\sigma)$ has coefficients in the integer ring of $E$ and they do not depend on $\mathcal{P}$. Fix an embedding $E \rightarrow \CC$. The absolute values of the complex roots of the characteristic polynomial of $\sigma$ are equal to ${\rm N}({\mathcal{L}})^{ 1/2 }=\lambda^{ 1/2 }$. 
\end{cor}
\DIM The first part of the previous corollary follows from \cite[2.1.2]{Rib}. Observe that the roots of the characteristic polynomial $\rho_{\mathcal{P}} (\sigma)$ are also roots of the characteristic polynomial of $\rho_p(\sigma)$, and so the bound for the absolute values of the roots is a consequence of Weil's results  (cf. Theorem \ref{WeilBound}). See also \cite[11.9 and 11.10]{Shimura} for further details. 
\CVD

Now we can prove the Corollary \ref{BoundIsogeny}.

\begin{cor}[{\rm Corollary \ref{BoundIsogeny}}]
Let $k$ and $E$ be number fields, denote $d=[E:\Q]$. Let $\A/k$ be an abelian variety of dimension $d$ such that $E$ embeds in $ {\rm End}_k(\A) \otimes \Q$. Let $p$ be a prime number. As before, denote $C = 2^{( d+1 )^{2 ln ( d ) / ln ( 2 )} + 1}$. Suppose that : 
\begin{enumerate}
\item  $k \cap \Q ( \zeta_p ) = \Q$;
\item $p$ is a good prime for $\A$;
\item $p\geq 3d+1$;
\item For each place $\P$ of $E$ dividing $p$ the determinant of the representation $\Gal (\overline{k}/k) \rightarrow \GL_{2}(E_{\mathcal{P}})$ is given by the cyclotomic character $\chi_{p}$;
\item $p$ is not ramified in $E$;
\item There exists a place $\mathcal{L}$ of $k$, not dividing $p$ such that $\A$ admits good reduction in $\mathcal{L}$;
\item Denote $\lambda$ the absolute norm of ${\mathcal{L}}$, that is the cardinal of the residue field of $\mathcal{L}$. For each place $\mathcal{P}$ of $E$ dividing $p$, the absolute norm $ {\rm N }(\P)$ of $\P$ satisfies the inequality ${\rm N }(\P) > (1+\lambda^{C/2})^{2d}$. 
\end{enumerate}
Then for all $n \in \N$ we have $H^1_{{\rm loc}} ( \Gal ( k ( \A[p^n] ) / k ) , \A[p^n] ) ) = 0$. Therefore, for all $n \in \N$, the local-global divisibility by $p^n$ holds for $\A ( k )$. If moreover $\A$ is principally polarized, $\Sh ( \A / k )$ is $p$-divisible in $H^1 ( k , \A )$.
\end{cor}

\DIM Suppose that $\mathcal{A}$, $E$, $\mathcal{L}$ and $p$ satisfy all the conditions 1-6 of Corollary \ref{BoundIsogeny} and that $H^1_{\rm loc} ( \Gal ( k ( \A[p^n] ) / k ) , \A[p^n] ) ) \neq 0$. By Lemma \ref{lem22}, there exists a prime $\mathcal{P}$ of $E$ dividing $p$ such that $H^1_{{\rm loc}} ( G_{\P , n} , \A[\P^n] ) \neq 0$. We will show that $N(P) \leq (1+\lambda^{C/2})^{2d}$.   
As $p$ is a good prime for $\mathcal{A}$, we have the representation 
$$ \rho_{\mathcal{P}} \colon \Gal(\overline{k} /k) \rightarrow \GL_2(\O_{E, \mathcal{P}}) $$
and  $\overline{\rho}_{\mathcal{P}}$, the reduction of $\rho_{\mathcal{P}}$ modulo $\P \O_{E , \P}$, is the representation given by the Galois action on the $\mathcal{P}$-torsion points of $\mathcal{A}$: 
$$ \overline{\rho}_{\mathcal{P}} \colon \Gal(\overline{k} /k) \rightarrow \GL_2(\O_{E, \mathcal{P}} /\mathcal{P}\O_{E, \mathcal{P}}  ). $$
By Proposition \ref{pro21} and Corollary \ref{cor36} there exists an isomorphism over $\F_q$ between $\GL_2(\mathcal{A}[\mathcal{P}])$ and $\GL_2(\F_q)$ such that any element $\sigma$ of the image of $\overline{\rho}_{\mathcal{P}} $ is triangular superior with diagonal coefficients $\alpha(\sigma)$, $\beta(\sigma)$. Moreover either $\alpha$ or $\beta$ has order $c$ for a constant $c\le C$ (let us suppose that it is $\alpha$). As the determinant of $\overline{\rho}_\P$ is $\overline{\chi}_p$, the diagonal coefficients are $\alpha $ and $\alpha^{-1} \overline{\chi}_p$. For $\sigma\in \Gal ( k ( \A[p] ) / k )$ the trace of $\overline{\rho}_{\P} (\sigma)$ is given by:
\begin{equation}
\label{PidFlo01}
\alpha(\sigma)+\alpha^{-1}(\sigma)\overline{\chi}_p(\sigma).
\end{equation}
As before fix a place $\overline{\mathcal{L}}$ of $\overline{\Q}$ over $\mathcal{L}$. Let $\sigma \in D(\overline{\mathcal{L}})$ such that the reduction of $\sigma$ in $\Gal ( \overline{\F}_{\lambda} / \F_{\lambda}) $ is the Frobenius morphism. We denote also by $\sigma$ its restriction in $\Gal ( k ( \A[p] ) / k )$.
Let $a_{\mathcal{L}}$ be the trace of  $\rho_{\mathcal{P}} (\sigma^{c})$. By Corollary \ref{Weil}, we have that the characteristic polynomial of $\rho_{\mathcal{P}} (\sigma)$ has coefficients in $O_E$. As a consequence it is also true for the characteristic polynomial of $\rho_{\mathcal{P}} (\sigma^{c})$. So $a_{\mathcal{L}}$ is an element of $\O_E$. The previous formula (\ref{PidFlo01}) allows us to compute $a_{\mathcal{L}}$ modulo $\mathcal{P} $. 
$$a_{\mathcal{L}} \equiv \alpha(\sigma^c)+\alpha^{-1}(\sigma^c)\overline{\chi}_p(\sigma^c)   \mod \mathcal{P}. $$
Observe that $\alpha^c=1$ and $\mathcal{L}$ does not divide $p$ so $\chi_{p} (\sigma)= \lambda$. We deduce that: 
$$ a_{\mathcal{L}} \equiv 1+\lambda^c \mod \mathcal{P}. $$
So 
$$ \P \ {\rm divides } \ (1+\lambda^c -a_{\mathcal{L}}). $$
In particular, if we denote by ${\rm N}_{E/\Q}$ the relative norm of $E$ over $\Q$, we have:
\begin{equation}
\label{PidFlo02}
{\rm N}(\P) \ {\rm divides } \  {\rm N}_{E/\Q}(1+\lambda^c -a_{\mathcal{L}} ).
\end{equation}
Now let $\tau$ be an embedding $\tau \colon E\rightarrow \CC$. For Corollary \ref{Weil} all the complex roots of the characteristic polynomial $\rho_{\mathcal{P}} (\sigma)$ have absolute value equal to $\lambda^{1/2}$, independently on the choice of an embedding $E\rightarrow \CC$. We deduce that all the complex roots of the characteristic polynomial $\rho_\P (\sigma^c)$ have absolute value equal to $\lambda^{c/2}$. So we have for any embedding $\tau \colon E\rightarrow \CC$ : 
$$ \vert \tau(a_{\mathcal{L}})  \vert \le 2\lambda^{c/2 },$$
that gives us an upper bound: 
$$ 0< \vert \tau (1+\lambda^c -a_{\mathcal{L}}) \vert \le (1+\lambda^{c/2})^2. $$
So we deduce (recall that $[E:\Q]=d$):
\begin{equation}
\label{eqWeil}
0< \vert {\rm N}_{E/\Q}(1+\lambda^c -a_{\mathcal{L}} ) \vert \le (1+\lambda^{c/2})^{2d}.
\end{equation}
By combining this inequality (\ref{eqWeil}), and the formula (\ref{PidFlo02}) using the fact that $c\le C$, we deduce 
$$ {\rm N}(\P) \leq (1+\lambda^{C/2})^{2d}.
 $$
\CVD

\end{document}